\documentclass{amsart}
\issueinfo{00}
{0} 
{0} 
{2007} 
\usepackage{xy}
\xyoption{all}
\usepackage{graphicx}
\usepackage{amsmath}
\newtheorem{theorem}{Theorem}[section]

\newtheorem{corollary}[theorem]{Corollary}

\newtheorem{lemma}[theorem]{Lemma}
\newtheorem{proposition}[theorem]{Proposition}
\theoremstyle{remark}
\newtheorem{remark}[theorem]{Remark}
\numberwithin{equation}{section}

\begin{document}
\title[Some geometric structures associated to a $k$-symplectic manifold]{Some geometric structures associated to a $k$-symplectic manifold}

\author[A. M. Blaga]{Adara M. Blaga}
 \address{Facultatea de Matematic\u{a} \c{s}i Informatic\u{a},
Universitatea de Vest din Timi\c{s}oara, Bd. Vasile P\^{a}rvan 4,
300223 Timi\c{s}oara, Romania}
 \email{adara@math.uvt.ro}

\author[B. Cappelletti Montano]{Beniamino Cappelletti Montano}
 \address{Dipartimento di Matematica,
Universit\`{a} degli Studi di Bari, Via Edoardo Orabona 4, 70125
Bari, Italy}
 \email{cappelletti@dm.uniba.it}

\subjclass[2000]{Primary 53C12, 53C15, Secondary 53B05, 53D05,
57R30.}

\keywords{$k$-symplectic structures, Ehresmann connections,
Lagrangian foliations, characteristic classes}

\begin{abstract}
A canonical connection is attached to any $k$-symplectic manifold.
We study the properties of this connection and its geometric
applications to $k$-symplectic manifolds. In particular, we prove
that, under some natural assumptions, any $k$-symplectic manifold
admits an Ehresmann connection, discuss some corollaries of this
result and find vanishing theorems for characteristic classes on a
$k$-symplectic manifold.
\end{abstract}

\maketitle

\section{Introduction}
The theory of $k$-symplectic manifolds was initiated by A. Awane
(\cite{awane1}), who defined a $k$-symplectic structure on an
$n(k+1)$-dimensional smooth manifold $M$ as an $n$-codimensional
foliation $\mathcal F$ and a system of $k$ closed $2$-forms
vanishing on the subbundle of \ $TM$ \ defined by $\mathcal F$ \
with transversal characteristic spaces \ (for a precise definition
see $\S$ \ref{sezionepreliminari}). The study of these structures
was motivated by some mathematical and physical considerations, like
the local study of Pfaffian systems and Nambu's statistical
mechanics. But the interest on $k$-symplectic geometry has increased
especially in recent years by the awareness of its relationship with
polysymplectic (or multisymplectic) and $n$-symplectic geometry, and
their applications in field theory (cf. \cite{deleon1},
\cite{merino}, \cite{rey}). In fact the $k$-symplectic formalism is
the generalization to field theories of the standard symplectic
formalism in mechanics, which is the geometric framework for
describing most of autonomous mechanical systems. Especially it can
be used for giving a geometric description of first order field
theories in which the  Lagrangian and Hamiltonian depend on the
first jet (prolongation) of the field.

The definition of a $k$-symplectic manifold is a  generalization of
the notion of a symplectic manifold foliated by a Lagrangian
foliation. Thus it is a natural question whether one can define an
appropriate analogue of the well-known notion of bi-Lagrangian
structure to the context of $k$-symplectic geometry. We recall that
an almost bi-Lagrangian manifold is a symplectic manifold
$(M^{2n},\omega)$ endowed with two transversal Lagrangian
distributions $L_1$ and $L_2$. When $L_1$ and $L_2$ are both
integrable, we speak of bi-Lagrangian manifold. The peculiarity of
these geometric structures is that a canonical symplectic connection
can be attached to them. This connection was introduced by H. Hess
(\cite{hess}), who was working in geometric quantization, and later
on its important geometric properties were pointed out by N. B.
Boyom (\cite{boyom}) and I. Vaisman (\cite{vaisman1},
\cite{vaisman2}).

In this work we consider the $k$-symplectic analogue of
bi-Lagrangian structure and attach to a such $k$-symplectic manifold
a canonical connection which plays the same role in $k$-symplectic
geometry as the Hess connection. Moreover we define on
 a $k$-symplectic manifold a family of tensor fields which can be
 thought as the proper generalization in this setting of almost
 K\"{a}hler structures, and we prove that under some integrability
 assumptions, the above connection coincides with the Levi-Civita
 connection of a suitable compatible metric. Finally, as an
 application, we prove that under some certain natural assumption, any
 $k$-symplectic manifold admits an Ehresmann connection and we
 deduce some geometric and topological properties on the
 $k$-symplectic manifold in question.

\section{$k$-symplectic structures}\label{sezionepreliminari}

A \emph{$k$-symplectic manifold} (cf. \cite{awane1}, \cite{puta}) is
a smooth manifold $M$ together with $k$ closed $2$-forms
$\omega_1,\ldots,\omega_k$ such that
\begin{enumerate}
    \item $C_x\left(\omega_1\right)\cap\cdots\cap
    C_x\left(\omega_k\right)=\left\{0\right\}$,
    \item $\omega_\alpha\left(X,X'\right)=0$ for any
    $X,X'\in\Gamma\left(T\mathcal{F}\right)$ and for all $\alpha\in\left\{1,\ldots,k\right\}$,
\end{enumerate}
where $C_x\left(\omega\right)=\left\{v\in T_xM :
\omega_x\left(v,w\right)=0 \textrm{ for any } w\in T_xM\right\}$ and
$\mathcal{F}$ is an $nk$-dimensional foliation on $M$. It follows
that $\dim\left(M\right)=n\left(k+1\right)$. We will usually denote
by $L$ the tangent bundle of the foliation $\mathcal{F}$. In terms
of $G$-structures, a $k$-symplectic manifold can be defined by an
integrable $\emph{{Sp}}(k,n;\mathbb{R})$-structure, where
$\emph{{Sp}}(k,n;\mathbb{R})$ denotes the $k$-symplectic group,
defined by the set of matrices of the following type
\begin{equation*}
\left(
  \begin{array}{cccc}
    T &  & 0 & S_1 \\
     & \ddots &  & \vdots \\
     &  & T & S_k \\
    0 &  &  & ^t(T^{-1}) \\
  \end{array}
\right)
\end{equation*}
where $T\in \emph{Gl}(n;\mathbb{R})$ and $S_1, \ldots, S_k$ are
$n\times n$ real matrices such that $T ^t S_\alpha = S_\alpha ^t T$
for all $\alpha\in\left\{1,\ldots,k\right\}$. The canonical model of
these structures is the $k$-cotangent bundle $(T^1_k)^\ast N$ of an
arbitrary manifold $N$, which can be identified with the vector
bundle $J^1(N,\mathbb{R}^k)_0$ whose total space is the manifold of
$1$-jets of maps with target $0\in\mathbb{R}^k$, and projection
$\tau^\ast(j^1_{x,0}\sigma)=x$. In this case, identifying
$(T^1_k)^\ast N$ with the Whitney sum of $k$ copies of $T^\ast N$,
$(T^1_k)^\ast N\cong T^\ast N\oplus\cdots\oplus T^\ast N$,
$j_{x,0}\sigma\mapsto(j^{1}_{x,0}\sigma^1,\ldots,j^{k}_{x,0}\sigma^k)$,
where
$\sigma^\alpha=\pi_\alpha\circ\sigma:N\longrightarrow\mathbb{R}$ is
the $\alpha$-th component of $\sigma$, the $k$-symplectic structure
on $(T^1_k)^\ast N$  is given by
$\omega_\alpha=(\tau_\alpha^\ast)^\ast(\omega_0)$ and ${T\mathcal
F}_{j_{x,0}^1\sigma}=\ker(\tau^\ast)_\ast(j_{x,0}^1\sigma)$, where
$\tau_\alpha^\ast:(T^1_k)^\ast N\longrightarrow T^\ast N$ is the
projection on the $\alpha$-th copy $T^\ast N$ of $(T^1_k)^\ast N$
and $\omega_0$ is the standard symplectic structure on $T^\ast N$.

Returning to the general case of an arbitrary $k$-symplectic
manifold $(M,\omega_\alpha,{\mathcal F})$, for each
$\alpha\in\left\{1,\ldots,k\right\}$ we set
\begin{equation}\label{n-lagrangian}
L_{{\alpha}_x}:=\bigcap_{\beta\neq\alpha}C_{x}\left(\omega_\beta\right).
\end{equation}
Then we have (\cite{awane2}):
\begin{enumerate}
    \item[(a)] for each $\alpha\in\left\{1,\ldots,k\right\}$ the distribution
    $L_\alpha=(L_{{\alpha}_x})_{x\in M}$ is integrable (we denote by $\mathcal{F}_\alpha$ the
    foliation integral to $L_\alpha$);
    \item[(b)] $L=L_1\oplus\cdots\oplus L_k$;
    \item[(c)] for each $\alpha\in\left\{1,\ldots,k\right\}$ the map
    $i_\alpha:L_\alpha\longrightarrow ({N\mathcal F})^\ast$, $X\mapsto
    i_X\omega_\alpha$, is an isomorphism, where $N\mathcal F$ denotes the
    normal bundle of $\mathcal F$.
\end{enumerate}

The standard Darboux theorem for Lagrangian foliations holds also
for $k$-symplectic manifolds:

\begin{theorem}[\cite{awane1}]\label{darboux}
About any point of a $k$-symplectic manifold
$\left(M,\omega_\alpha,\mathcal F\right)$, $\alpha \in \{1,...,k\}$,
there exist local coordinates
$\left\{x_1,\ldots,x_n,y_{1},\ldots,y_{kn}\right\}$ such that
$\omega_\alpha=\sum_{i=1}^{n}dx_i\wedge dy_{(\alpha-1)n+i}$ and
$\mathcal F$ is described by the equations
$\left\{x_i=\textrm{const.}\right\}$. In particular, for each
$\alpha\in\left\{1,\ldots,k\right\}$, $L_\alpha$ is generated by
$\frac{\displaystyle
\partial}{\displaystyle \partial y_{(\alpha-1)n+1}},\ldots,\frac{\displaystyle \partial}{\displaystyle \partial
y_{\alpha n}}$.
\end{theorem}

Recall that a vector field $X$ on a symplectic manifold
$(M^{2n},\omega)$ is said to be symplectic if ${\mathcal
L}_{X}\omega=0$. For $k$-symplectic manifolds we prove the following
lemma which will be useful in the sequel.

\begin{lemma}\label{liederivative}
In any $k$-symplectic manifold, ${\mathcal L}_{X}\omega_\alpha=0$,
for any $X\in\Gamma\left(L_\beta\right)$ with $\alpha\neq\beta$.
\end{lemma}
\begin{proof}
Using the Cartan formula for the Lie derivative, we have ${\mathcal
L}_{X}\omega_\alpha=i_{X}d\omega_\alpha+di_{X}\omega_\alpha=
di_{X}\omega_\alpha$, since $\omega_\alpha$ is closed. But, for any
$V\in\Gamma\left(TM\right)$,
$i_{X}\omega_\alpha\left(V\right)=2\omega_{\alpha}\left(X,V\right)=0$
from the definition of $L_\alpha$.
\end{proof}

\section{A canonical connection on $k$-symplectic manifolds}\label{canonicalconnection}

Let $\left(M,\omega_\alpha,\mathcal F\right)$,
$\alpha\in\left\{1,\ldots,k\right\}$, be a $k$-symplectic manifold.
In what follows, $Q$ will denote an $n$-dimensional integrable
distribution on $M$ transversal to $\mathcal{F}$ such that
\begin{enumerate}
    \item[(i)] $\omega_\alpha\left(Y,Y'\right)=0$ for any
    $Y,Y'\in\Gamma\left(Q\right)$ and for all
    $\alpha\in\left\{1,\ldots,k\right\}$,
    \item[(ii)]  $\left[X,Y\right]\in\Gamma\left(L_\alpha\oplus Q\right)$
    for any $X\in\Gamma\left(L_\alpha\right)$ and for any $Y\in\Gamma\left(Q\right)$.
\end{enumerate}
Occasionally, we will denote by $\mathcal G$ the foliation integral
to $Q$.

The geometric interpretation of the condition (i) is that, for each
$\alpha\in\left\{1,\ldots,k\right\}$ and for any $x\in M$, $Q_x$ is
a Lagrangian subspace of the symplectic vector space
$(L_{{\alpha}_{x}}\oplus Q_{x},\omega_{{\alpha}_{x}})$. The
condition (ii) is more technical; it will be essential for proving
some preliminary results, like the following Lemma \ref{lemmadue},
and then for the generalization of the Hess's construction to the
$k$-symplectic setting. Its geometric meaning is that for each fixed
$\alpha\in\left\{1,\ldots,k\right\}$, the subbundle $L_\alpha\oplus
Q$ is integrable, hence it defines a foliation whose leaves are
symplectic manifolds with respect to the restriction of the
$k$-symplectic form $\omega_\alpha$ to the leaves. We also have that
$(L_\alpha,Q)$ is a bi-Lagrangian structure on the leaves of the
foliation defined by $L_\alpha\oplus Q$.

A simple example of a $k$-symplectic manifold endowed with a
transversal integrable distribution verifying (i) and (ii) is given
by $\mathbb{R}^{n(k+1)}$ with its standard $k$-symplectic structure
given by Theorem \ref{darboux} and taking as $Q$ the distribution
spanned by $\frac{\displaystyle \partial}{\displaystyle \partial
x_1},\ldots,\frac{\displaystyle \partial}{\displaystyle \partial
x_n}$.

We also remark that the splitting $TM=L\oplus
Q=L_1\oplus\cdots\oplus L_k\oplus Q$ induces a canonical isomorphism
between $Q$ and $N\mathcal F:=TM/L$, the normal bundle to the
foliation $\mathcal F$. In particular, it follows that
$Q^{\ast}=\textrm{ann}(L)$ and, arguing in the same way for the
foliation $\bigoplus_{\beta\neq\alpha}L_\beta\oplus Q$, we get that
$L_\alpha^\ast=\textrm{ann}(\bigoplus_{\beta\neq\alpha}L_\beta\oplus
Q)$, for each $\alpha\in\left\{1,\ldots,k\right\}$. Taking into
account these remarks, we can prove the following preliminary
lemmas:

\begin{lemma}\label{lemmauno}
Let $X,X'\in\Gamma\left(L\right)$. For each
$\alpha\in\left\{1,\ldots,k\right\}$, the map
\begin{equation*}
\varphi^{XX'}_\alpha:
V\mapsto\left({\mathcal{L}}_{X}i_{X'}\omega_\alpha\right)\left(V\right)=
X\left(\omega_\alpha\left(X',V\right)\right)-\omega_\alpha\left(X',\left[X,V\right]\right),
\end{equation*}
for any $V\in \Gamma(TM)$, belongs to $Q^\ast$.
\end{lemma}
\begin{proof}
For any $X''\in\Gamma\left(L\right)$,
$\left({\mathcal{L}}_{X}i_{X'}\omega_\alpha\right)\left(X''\right)=X\left(\omega_\alpha\left(X',X''\right)\right)-\omega_\alpha\left(X',\left[X,X''\right]\right)=0$,
from which, since $Q^\ast=\textrm{ann}\left(L\right)$, we get the
result.
\end{proof}

\begin{lemma}\label{lemmadue}
Let $Y,Y'\in\Gamma\left(Q\right)$. For each
$\alpha\in\left\{1,\ldots,k\right\}$, the map
\begin{equation*}
\psi^{YY'}_\alpha:
V\mapsto\left({\mathcal{L}}_{Y}i_{Y'}\omega_\alpha\right)\left(V\right)=
Y\left(\omega_\alpha\left(Y',V\right)\right)-\omega_\alpha\left(Y',\left[Y,V\right]\right),
\end{equation*}
for any $V\in \Gamma(TM)$, belongs to $L_\alpha^\ast$.
\end{lemma}
\begin{proof}
Since
$L_\alpha^\ast=\textrm{ann}(\bigoplus_{\beta\neq\alpha}L_\beta\oplus
Q)$, we have to prove that
$\left({\mathcal{L}}_{Y}i_{Y'}\omega_\alpha\right)\left(X\right)=0$
and
$\left({\mathcal{L}}_{Y}i_{Y'}\omega_\alpha\right)\left(Y''\right)=0$
for any $X\in\Gamma\left(L_\beta\right)$, $\beta\neq\alpha$, and for
any $Y''\in\Gamma\left(Q\right)$. Indeed,
$\left({\mathcal{L}}_{Y}i_{Y'}\omega_\alpha\right)\left(X\right)=Y\left(\omega_\alpha\left(Y',X\right)\right)-\omega_\alpha\left(Y',\left[Y,X\right]\right)=0$
by the definition of $L_\beta$ and by (ii). Next,
$\left({\mathcal{L}}_{Y}i_{Y'}\omega_\alpha\right)\left(Y''\right)=Y\left(\omega_\alpha\left(Y',Y''\right)\right)-\omega_\alpha\left(Y',\left[Y,Y''\right]\right)=0$
by (i) and by the integrability of $Q$.
\end{proof}

\begin{theorem}\label{connection}
Let $\left(M,\omega_\alpha,\mathcal F\right)$, $\alpha \in
\{1,...,k\}$, be a $k$-symplectic manifold and let $Q$ be an
integrable distribution supplementary to $T\mathcal F$ verifying the
above conditions (i), (ii) and such that
$(i^{\ast}_{1})^{-1}(\psi^{YY'}_1)=\cdots=(i^{\ast}_{k})^{-1}(\psi^{YY'}_k)$
for any $Y,Y'\in\Gamma(Q)$, where $\psi^{YY'}_1,\ldots,\psi^{YY'}_k$
are the maps defined in Lemma \ref{lemmadue}. Then there exists a
unique connection $\nabla$ on $M$ satisfying the following
properties:
\begin{enumerate}
    \item  $\nabla{\mathcal F_\alpha}\subset{\mathcal F_\alpha}$ for each $\alpha\in\left\{1,\ldots,k\right\}$, and $\nabla Q\subset Q$,
    \item $\nabla\omega_1=\cdots=\nabla\omega_k=0$,
    \item $T\left(X,Y\right)=0$ for any $X\in\Gamma\left(L\right)$ and
    for any $Y\in\Gamma\left(Q\right)$,
\end{enumerate}
where $T$ denotes the torsion tensor field of $\nabla$.
\end{theorem}
\begin{proof}
According to the decomposition $TM=L_1\oplus\cdots\oplus L_k\oplus
Q$, we define a connection $\nabla^{L_\alpha}$ on each subbundle
$L_\alpha$, a connection $\nabla^Q$ on $Q$ and then we take the sum
of these connections for defining a global connection on $M$. Fix an
$\alpha\in\left\{1,\ldots,k\right\}$. We define
$\nabla^{L_\alpha}_{Y}X:=\left[Y,X\right]_{L_\alpha}$  for any
$X\in\Gamma\left(L_\alpha\right)$ and $Y\in\Gamma\left(Q\right)$.
Now we have to define $\nabla^{L_\alpha}_{X}X'$ for
$X\in\Gamma\left(L\right)$, $X'\in\Gamma\left(L_\alpha\right)$.
Since $i_\alpha:L_\alpha\longrightarrow Q^\ast$ is an isomorphism
for any fixed $X\in \Gamma\left(L\right)$, $X'\in
\Gamma\left(L_\alpha\right)$, by Lemma \ref{lemmauno}, there exists
a unique section
$H_\alpha\left(X,X'\right)\in\Gamma\left(L_\alpha\right)$ such that
$i_\alpha\left(H_\alpha\left(X,X'\right)\right)=\varphi^{XX'}_\alpha$,
that is
$\omega_\alpha\left(H_\alpha\left(X,X'\right),Y\right)=X\left(\omega_\alpha\left(X',Y\right)\right)-\omega_\alpha\left(X',\left[X,Y\right]\right)$
for any $Y\in\Gamma\left(Q\right)$. We set
$\nabla^{L_\alpha}_{X}X':=H_\alpha\left(X,X'\right)\in\Gamma\left(L_\alpha\right)$.
Now we define the connection $\nabla^Q$. For any
$X\in\Gamma\left(L\right)$ and $Y\in\Gamma\left(Q\right)$ we put
$\nabla^{Q}_{X}Y:=\left[X,Y\right]_Q$. It remains to define
$\nabla^Q_{Y}Y'$ for $Y,Y'\in\Gamma\left(Q\right)$. The isomorphism
$i_\alpha:L_\alpha\longrightarrow Q^\ast$ determines an isomorphism
 $i^{\ast}_\alpha$ between $Q$ and $L_\alpha^\ast$ such
that
$i^{\ast}_\alpha\left(Y\right)\left(X\right)=\omega_\alpha\left(Y,X\right)$.
Then, for any fixed $Y,Y'\in\Gamma\left(Q\right)$, by Lemma
\ref{lemmadue}, there exists a unique section
$H_\alpha\left(Y,Y'\right)\in \Gamma(Q)$ such that
$i^{\ast}_\alpha\left(H\left(Y,Y'\right)\right)=\psi^{YY'}_\alpha$,
that is
$\omega_\alpha\left(H_\alpha\left(Y,Y'\right),X\right)=Y\left(\omega_\alpha\left(Y',X\right)\right)-\omega_\alpha\left(Y',\left[Y,X\right]\right)$
for any $X\in\Gamma\left(L_\alpha\right)$. Moreover, our assumption
ensures that
$H_1\left(Y,Y'\right)=\cdots=H_k\left(Y,Y'\right)=:H\left(Y,Y'\right)$.
We set $\nabla^Q_{Y}Y':=H\left(Y,Y'\right)\in\Gamma\left(Q\right)$.
Now we prove that $\nabla^Q$ is a connection on $Q$ and, for each
$\alpha\in\left\{1,\ldots,k\right\}$, $\nabla^{L_\alpha}$ is a
connection on $L_\alpha$. For any $X\in\Gamma\left(L\right)$,
$Y\in\Gamma\left(Q\right)$ and $f\in C^\infty\left(M\right)$ we have
\begin{gather*}
\nabla^Q_{fX}Y=\left[fX,Y\right]_Q=f\left[X,Y\right]_Q-Y\left(f\right)X_Q=f\left[X,Y\right]_Q=f\nabla^Q_XY,\\
\nabla^Q_{X}\left(fY\right)=\left[X,fY\right]_Q=f\left[X,Y\right]_Q+X\left(f\right)Y_Q=f\nabla^Q_XY+X\left(f\right)Y,
\end{gather*}
and, for any $X\in\Gamma\left(L_\alpha\right)$,
$Y,Y'\in\Gamma\left(Q\right)$,
\begin{align*}
\omega_\alpha(\nabla^Q_{fY}Y',X)&=\omega_\alpha\left(H_\alpha\left(fY,Y'\right),X\right)\\
&=fY\left(\omega_\alpha\left(Y',X\right)\right)-\omega_\alpha\left(Y',\left[fY,X\right]\right)\\
&=fY\left(\omega_\alpha\left(Y',X\right)\right)-f\omega_\alpha\left(Y',\left[Y,X\right]\right)+X\left(f\right)\omega_\alpha\left(Y',Y\right)\\
&=f\omega_\alpha\left(H_\alpha\left(Y,Y'\right),X\right)\\
&=\omega_\alpha(f\nabla^Q_YY',X),
\end{align*}
from which we get $\nabla^Q_{fY}Y'=f\nabla^Q_{Y}Y'$. Moreover,
\begin{align*}
\omega_\alpha(\nabla^Q_{Y}(fY'),X)&=\omega_\alpha\left(H\left(Y,fY'\right),X\right)\\
&=Y\left(\omega_\alpha\left(fY',X\right)\right)-\omega_\alpha\left(fY',\left[Y,X\right]\right)\\
&=fY\left(\omega_\alpha\left(Y',X\right)\right)+Y\left(f\right)\omega_\alpha\left(Y',X\right)
-f\omega_\alpha\left(Y',\left[Y,X\right]\right)\\&=
f\omega_\alpha\left(H(Y,Y'),X\right)+Y(f)\omega_\alpha\left(Y',X\right)
\\&=\omega_\alpha(f\nabla^Q_{Y}Y'+Y\left(f\right)Y',X),
\end{align*}
from which we obtain
$\nabla^Q_{Y}\left(fY'\right)=f\nabla^Q_{Y}Y'+Y\left(f\right)Y'$.
Now we prove that $\nabla^{L_\alpha}$ is a connection on the
subbundle $L_\alpha$, for each $\alpha\in\left\{1,\ldots,k\right\}$.
As before it is easy to show that
$\nabla^{L_\alpha}_{fY}X=f\nabla^{L_\alpha}_{Y}X$ and
$\nabla^{L_\alpha}_{Y}\left(fX\right)=f\nabla^{L_\alpha}_{Y}X+Y\left(f\right)X$
for any $X\in\Gamma\left(L_\alpha\right)$ and
$Y\in\Gamma\left(Q\right)$. Then for any $X\in\Gamma\left(L\right)$,
$X'\in\Gamma\left(L_\alpha\right)$ and any
$Y\in\Gamma\left(Q\right)$
\begin{align*}
\omega_\alpha(\nabla^{L_{\alpha}}_{fX}X',Y)&=\omega_\alpha\left(H_{\alpha}\left(fX,X'\right),Y\right)\\
&=fX\left(\omega_\alpha\left(X',Y\right)\right)-\omega_\alpha\left(X',\left[fX,Y\right]\right)\\
&=fX\left(\omega_\alpha\left(X',Y\right)\right)-f\omega_\alpha\left(X',\left[X,Y\right]\right)+Y\left(f\right)\omega_\alpha\left(X',X\right)\\
&=f\omega_\alpha\left(H_{\alpha}\left(X,X'\right),Y\right)\\
&=\omega_\alpha(f\nabla^{L_{\alpha}}_{X}X',Y),
\end{align*}
from which we get
$\nabla^{L_{\alpha}}_{fX}X'=f\nabla^{L_{\alpha}}_{X}X'$. Moreover,
\begin{align*}
\omega_\alpha(\nabla^{L_{\alpha}}_{X}\left(fX'\right),Y)&=\omega_\alpha\left(H_{\alpha}\left(X,fX'\right),Y\right)\\
&=X\left(\omega_\alpha\left(fX',Y\right)\right)-\omega_\alpha\left(fX',\left[X,Y\right]\right)\\
&=fX\left(\omega_\alpha\left(X',Y\right)\right)+X\left(f\right)\omega_\alpha\left(X',Y\right)-f\omega_\alpha\left(X',\left[X,Y\right]\right)\\
&=f\omega_\alpha\left(H_{\alpha}\left(X,X'\right),Y\right)+X\left(f\right)\omega_\alpha\left(X',Y\right)\\
&=\omega_\alpha(f\nabla^{L_{\alpha}}_{X}X'+X\left(f\right)X',Y)
\end{align*}
from which we get
$\nabla^{L_{\alpha}}_{X}\left(fX'\right)=f\nabla^{L_{\alpha}}_{X}X'+X\left(f\right)X'$.
Therefore we can define a global connection on $M$ putting, for
any $V,W\in\Gamma\left(TM\right)$,
\begin{equation}
\nabla_{V}W=\nabla^{L_1}_{V}W_{L_1}+\cdots+\nabla^{L_k}_{V}W_{L_k}+\nabla^{Q}_{V}W_{Q}.
\end{equation}
Now we prove that the connection $\nabla$ satisfies (1)--(3). By
construction $\nabla$ preserves the distributions $L_\alpha$ and
$Q$. Then, by (1) we have that, obviously,
$\left(\nabla_{V}\omega_\alpha\right)\left(X,X'\right)=0$ for any
$X,X'\in\Gamma\left(L\right)$ and $V\in\Gamma\left(TM\right)$. For
the same reason,
$\left(\nabla_{V}\omega_\alpha\right)\left(Y,Y'\right)=0$ for any
$Y,Y'\in\Gamma\left(Q\right)$ and $V\in \Gamma(TM)$. Now, let
$X\in\Gamma\left(L\right)$, $X'\in\Gamma\left(L_\alpha\right)$ and
$Y\in\Gamma\left(Q\right)$. Then
\begin{align*}
\left(\nabla_{X}\omega_\alpha\right)\left(X',Y\right)&=X\left(\omega_\alpha\left(X',Y\right)\right)-\omega_\alpha\left(H\left(X,X'\right),Y\right)-\omega_\alpha(X',\left[X,Y\right]_Q)\\
&=X\left(\omega_\alpha\left(X',Y\right)\right)-X\left(\omega_\alpha\left(X',Y\right)\right)+\omega_\alpha\left(X',\left[X,Y\right]\right)\\
&\quad-\omega_\alpha\left(X',\left[X,Y\right]\right)=0.
\end{align*}
Moreover, for any $\beta\neq\alpha$
$\left(\nabla_{X}\omega_\beta\right)\left(X',Y\right)=0$ because
$\nabla_{X}X'\in\Gamma\left(L_\alpha\right)$. Finally, for any
$X'\in\Gamma\left(L_\alpha\right)$ and
$Y,Y'\in\Gamma\left(Q\right)$,
\begin{align*}
\left(\nabla_{Y}\omega_\alpha\right)\left(X',Y'\right)&=Y\left(\omega_\alpha\left(X',Y'\right)\right)-\omega_\alpha(\left[Y,X'\right]_{L_\alpha},Y')-\omega_\alpha\left(X',H\left(Y,Y'\right)\right)\\
&=Y\left(\omega_\alpha\left(X',Y'\right)\right)-\omega_\alpha\left(\left[Y,X'\right]_{L_{\alpha}},Y'\right)+Y\left(\omega_\alpha\left(Y',X'\right)\right)\\
&\quad-\omega_\alpha\left(Y',\left[Y,X'\right]\right)=0.
\end{align*}
Thus we conclude that
$\left(\nabla_{V}\omega_\alpha\right)\left(X,Y\right)=0$ for any
$X\in\Gamma\left(L_\alpha\right)$, $Y\in\Gamma\left(Q\right)$ and
$V\in\Gamma\left(TM\right)$. Analogously, one can compute for all
the other cases, concluding that $\nabla\omega_\alpha=0$ for all
$\alpha\in\left\{1,\ldots,k\right\}$. Finally, for any
$X\in\Gamma\left(L_\alpha\right)$ and $Y\in\Gamma\left(Q\right)$ we
have
$T\left(X,Y\right)=\left[X,Y\right]_Q-\left[Y,X\right]_{L_{\alpha}}-\left[X,Y\right]=\left[X,Y\right]_{L_\alpha\oplus
Q}-\left[X,Y\right]=0$, since by (ii)
$\left[X,Y\right]\in\Gamma(L_{\alpha}\oplus Q)$. It remains to prove
the uniqueness of this connection up to the properties (1)--(3). Let
$X\in\Gamma\left(L\right)$ and $Y\in\Gamma\left(Q\right)$. For any
$X'\in\Gamma\left(L\right)$ we have, by (1) and (3),
$\omega_\alpha\left(\nabla_{X}Y,X'\right)=\omega_\alpha\left(\nabla_{Y}X+\left[X,Y\right],X'\right)=\omega_\alpha\left(\left[X,Y\right],X'\right)$,
for all $\alpha\in\left\{1,\ldots,k\right\}$, from which we get
$\nabla_{X}Y=\left[X,Y\right]_Q$. Then, using (3) again, we obtain
$\nabla_{Y}X=\left[Y,X\right]_L$. Moreover, for any
$X\in\Gamma\left(L\right)$, $X'\in\Gamma\left(L_\alpha\right)$ and
$Y\in\Gamma\left(Q\right)$ by (2) we have
$\omega_\alpha\left(\nabla_{X}X',Y\right)=X\left(\omega_\alpha\left(X',Y\right)\right)-\omega_\alpha\left(X',\nabla_{X}Y\right)=X\left(\omega_\alpha\left(X',Y\right)\right)-\omega_\alpha(X',\left[X,Y\right]_Q)=X\left(\omega_\alpha\left(X',Y\right)\right)-\omega_\alpha\left(X',\left[X,Y\right]\right)=\omega_\alpha\left(H_\alpha\left(X,X'\right),Y\right)$,
from which, since
$\nabla_{X}X',H_\alpha\left(X,X'\right)\in\Gamma\left(L_\alpha\right)$,
we get $\nabla_{X}X'=H_\alpha\left(X,X'\right)$. Similarly, one can
find that $\nabla_{Y}Y'=H\left(Y,Y'\right)$ for any
$Y,Y'\in\Gamma\left(Q\right)$.
\end{proof}

\begin{proposition}\label{torsion}
The connection $\nabla$ defined in Theorem \ref{connection} is
torsion free along the leaves of the foliations $\mathcal F$ and
$\mathcal G$.
\end{proposition}
\begin{proof}
Let $X\in\Gamma\left(L_\beta\right)$ and
$X'\in\Gamma\left(L_\alpha\right)$ and assume that
$\alpha\neq\beta$. We have
$T\left(X,X'\right)=H_\alpha\left(X,X'\right)-H_\beta\left(X',X\right)-\left[X,X'\right]\in\Gamma\left(L\right)$.
Then for any $Y\in\Gamma\left(Q\right)$
\begin{align*}
\omega_\alpha\left(T\left(X,X'\right),Y\right)&=\omega_\alpha\left(H_\alpha\left(X,X'\right)-\left[X,X'\right],Y\right)\\
&=X\left(\omega_\alpha\left(X',Y\right)\right)-\omega_\alpha\left(X',\left[X,Y\right]\right)-\omega_\alpha\left(\left[X,X'\right],Y\right)\\
&=3d\omega_\alpha\left(X,X',Y\right)=0
\end{align*}
since each $\omega_\alpha$ is closed. Analogously,
$\omega_\alpha\left(T\left(X,X'\right),Y\right)=0$. Moreover, for
each $\gamma\neq\alpha,\beta$
\begin{equation*}
\omega_\gamma\left(T\left(X,X'\right),Y\right)=-\omega_\gamma\left(\left[X,X'\right],Y\right)=3d\omega_\gamma\left(X,X',Y\right)=0.
\end{equation*}
Then $T\left(X,X'\right)\in C\left(\omega_1\right)\cap\cdots\cap
C\left(\omega_k\right)=\left\{0\right\}$. If
$X,X'\in\Gamma\left(L_\alpha\right)$, we have
$T\left(X,X'\right)=H_\alpha\left(X,X'\right)-H_\alpha\left(X',X\right)-\left[X,X'\right]\in\Gamma\left(L_\alpha\right)$
and
\begin{align*}
\omega_\alpha\left(T\left(X,X'\right),Y\right)&=X\left(\omega_\alpha\left(X',Y\right)\right)-\omega_\alpha\left(X',\left[X,Y\right]\right)-X'\left(\omega_\alpha\left(X,Y\right)\right)\\
&\quad+\omega_\alpha\left(X,\left[X',Y\right]\right)-\omega_\alpha\left(\left[X,X'\right],Y\right)\\
&=3d\omega_\alpha\left(X,X',Y\right)=0,
\end{align*}
hence $T\left(X,X'\right)=0$. Analogously, one can prove that
$T\left(Y,Y'\right)=0$ for any $Y,Y'\in\Gamma\left(Q\right)$.
\end{proof}

\begin{proposition}\label{curvature}
The curvature tensor field of the connection $\nabla$ defined in
Theorem \ref{connection} vanishes along the leaves of the foliations
${\mathcal F}$ and $\mathcal G$.
\end{proposition}
\begin{proof}
For any $X,X'\in\Gamma\left(L\right)$ and
$Y\in\Gamma\left(Q\right)$, using the integrability of $L$, we have
\begin{align*}
R_{X, X'}Y&=\nabla_{X}\left[X',Y\right]_Q-\nabla_{X'}\left[X,Y\right]_Q-\nabla_{\left[X,X'\right]}Y\\
&=\nabla_{X}\left[X',Y\right]_Q-\nabla_{X'}\left[X,Y\right]_Q-\left[\left[X,X'\right],Y\right]_Q=0
\end{align*}
by the Jacobi identity. Then, for any $X,X'\in\Gamma\left(L\right)$
and $X''\in\Gamma\left(L_\alpha\right)$ we have
\begin{equation}
R_{X,
X'}X''=H_\alpha\left(X,H_\alpha\left(X',X''\right)\right)-H_\alpha\left(X',H_\alpha\left(X,X''\right)\right)-H_\alpha\left(\left[X,X'\right],X''\right)
\end{equation}
Now, for any $Y\in\Gamma\left(Q\right)$
\begin{align*}
\omega_\alpha(H_\alpha(X,H_\alpha(X',X'')),Y)&=X(\omega_\alpha(H_\alpha(X',X''),Y))-\omega_\alpha(H_\alpha(X',X''),[X,Y])\\
&=X(\omega_\alpha(H(X',X''),Y))-\omega_\alpha(H(X',X''),[X,Y])\\
&=X(X'(\omega_\alpha(X'',Y)))-X(\omega_\alpha(X'',[X',Y]))\\
&\quad-X'(\omega_\alpha(X'',[X',[X,Y]])+\omega_\alpha(X'',[X',[X,Y]]),
\end{align*}
\begin{align*}
\omega_\alpha(H_\alpha (X',H_\alpha(X,X'')),Y)&=
X'(\omega_\alpha(H_\alpha(X,X''),Y))-\omega_\alpha(H_\alpha(X,X''),[X',Y])\\
&=X'(\omega_\alpha(H(X,X''),Y))-\omega_\alpha(H(X,X''),[X',Y])\\
&=X'(X(\omega_\alpha(X'',Y)))-X'(\omega_\alpha(X'',[X,Y]))\\
&\quad-X(\omega_\alpha(X'',[X,[X',Y]]))+\omega_\alpha(X'',[X,[X',Y]])
\end{align*}
and
\begin{equation*}
\omega_\alpha\left(H_\alpha\left(\left[X,X'\right],X''\right),Y\right)=\left[X,X'\right]\left(\omega_\alpha\left(X'',Y\right)\right)-\omega_\alpha\left(X'',\left[\left[X,X'\right],Y\right]\right).
\end{equation*}
Therefore
\begin{align*}
\omega_\alpha\left(R_{X, X'}X'',Y\right)&=\left[X,X'\right]\left(\omega_\alpha\left(X'',Y\right)\right)+\omega_\alpha\left(X'',\left[X',\left[X,Y\right]\right]\right)-\omega_\alpha\left(X'',\left[X,\left[X',Y\right]\right]\right)\\
&\quad-\left[X,X'\right]\left(\omega_\alpha\left(X'',Y\right)\right)+\omega_\alpha\left(X'',\left[\left[X,X'\right],Y\right]\right)\\
&=\omega_\alpha\left(X'',\left[\left[X,X'\right],Y\right]+\left[\left[X',Y\right],X\right]+\left[\left[Y,X\right],X'\right]\right)=0
\end{align*}
by the Jacobi identity. This shows that $R_{X, X'}=0$ for any
$X,X'\in\Gamma\left(L\right)$. In the same way, one can prove the
flatness along the leaves of the foliation defined by $Q$.
\end{proof}

\begin{corollary}
The leaves of the foliations $\mathcal F$ and $\mathcal G$ admit a
canonical flat affine structure.
\end{corollary}

Now we give an interpretation of the connection stated in Theorem
\ref{connection} in terms of some geometric structures which can be
attached to a $k$-symplectic manifold. So let
$(M,\omega_\alpha,\mathcal F)$, $\alpha \in \{1,...,k\}$, be a
$k$-symplectic manifold and let $Q$ be a distribution transversal to
$\mathcal F$ such that $\omega_\alpha(Y,Y')=0$ for any
$Y,Y'\in\Gamma(Q)$. Assume that $M$ admits a Riemannian metric $g$
such that the distributions $L_1, \ldots, L_k, Q$ are mutually
orthogonal. For each $\alpha \in \{1,...,k\}$, since $\omega_\alpha$
is non-degenerate on $L_\alpha\oplus Q$,  one can find a linear map
$A_{\alpha}:L_\alpha\oplus Q\longrightarrow L_\alpha\oplus Q$ such
that $\omega_{\alpha}(X,Y)=g(X,A_{\alpha}Y)$, for any $X,Y\in
\Gamma(L_\alpha\oplus Q)$. The operator $A_{\alpha}$, $\alpha \in
\{1,\ldots,k\}$, is skew-symmetric and $A_{\alpha}A_{\alpha}^*$,
$\alpha \in \{1,\ldots,k\}$, is symmetric and positive definite,
thus it diagonalizes with positive eigenvalues
$(\lambda_{\alpha})_i$, $i \in \{1,\ldots,2n\}$,
$A_{\alpha}A_{\alpha}^*=B_{\alpha}
\textrm{diag}\{(\lambda_{\alpha})_1,\ldots,(\lambda_{\alpha})_{2n}\}B_{\alpha}^{-1}$.
Set $\sqrt{A_{\alpha}A_{\alpha}^*}:=B_{\alpha}
\textrm{diag}\{\sqrt{(\lambda_{\alpha})_1},\ldots,\sqrt{(\lambda_{\alpha})_{2n}}\}B_{\alpha}^{-1}$
which is also symmetric and positive definite. Set
\begin{equation*}
J_{\alpha}:=\left\{%
\begin{array}{ll}
    (\sqrt{A_{\alpha}A_{\alpha}^*})^{-1}A_{\alpha}, & \hbox{on $L_\alpha\oplus Q$;} \\
    0, & \hbox{on $L_\beta$, $\beta\neq\alpha$.} \\
\end{array}%
\right.
\end{equation*}
Then $(J_1,\ldots,J_k)$ is a family of endomorphisms of the tangent
space satisfying
\begin{enumerate}
  \item[(i)]  $L_\alpha=\bigcap_{\beta\neq\alpha}\ker(J_\beta)$,
  \item[(ii)] $J_{\alpha}^2=-I$ on $L_\alpha\oplus Q$ and $J_\alpha
L_\alpha=Q$, $J_\alpha Q=L_\alpha$,
  \item[(iii)] $\omega_\alpha(X,Y)=g(X,J_\alpha Y)$ for any $X,Y\in\Gamma(TM)$.
\end{enumerate}
Note also that the Riemannian metric $g$ satisfies  $g(J_\alpha X,
J_\alpha Y)=g(X,Y)$ for each $\alpha\in\left\{1,\ldots,k\right\}$
and for any $X,Y\in\Gamma(TM)$. We call $(J_1,\ldots,J_k,g)$ a
\emph{compatible almost $k$-K\"{a}hler structure}. Now assume to be
under the assumptions of Theorem \ref{connection}. Note that for
each $\alpha\in\left\{1,\ldots,k\right\}$, the leaves of the
foliation defined by $L_\alpha\oplus Q$, endowed with the tensor
fields induced by $J_\alpha$, are almost K\"{a}hler manifolds. Then
we have that $[J_\alpha,J_\alpha]=0$ if and only if each leaf of the
foliation $L_\alpha\oplus Q$ is K\"{a}hlerian. When
$[J_\alpha,J_\alpha]=0$, for each
$\alpha\in\left\{1,\ldots,k\right\}$, that is the leaves of all the
foliations $L_\alpha\oplus Q$ are K\"{a}hler manifolds, we say that
$(M,\omega_\alpha,{\mathcal F},J_\alpha,g)$ is a
\emph{$k$-K\"{a}hler manifold}. Then we have the following result.

\begin{theorem}
Let $(M,\omega_\alpha,{\mathcal F},J_\alpha,g)$,
$\alpha\in\left\{1,\ldots,k\right\}$, be a $k$-K\"{a}hler manifold.
If the Levi-Civita connection $\nabla^g$ preserves the distributions
$L_\alpha$, then it preserves also $Q$ and it coincides with the
canonical connection $\nabla$.
\end{theorem}
\begin{proof}
We show that the Levi-Civita connection $\nabla^g$ satisfies the
properties (1), (2), (3) which, according to Theorem
\ref{connection}, define uniquely the canonical connection $\nabla$.
First of all we prove that $\nabla^g$ preserves $Q$. Let
$Y\in\Gamma\left(Q\right)$. Then, since $\nabla^g g=0$, for any
$V\in\Gamma\left(TM\right)$ and $X\in\Gamma\left(T\mathcal
F\right)$, we have
\begin{equation*}
0=(\nabla^{g}_{V}g)\left(X,Y\right)=V\left(g\left(X,Y\right)\right)-g(\nabla^{g}_{V}X,Y)-g(X,\nabla^{g}_{V}Y)=-g(X,\nabla^{g}_{V}Y),
\end{equation*}
since $\nabla^g{\mathcal F}\subset{\mathcal F}$. Thus
$\nabla^{g}Q\subset Q$. Finally we have to prove that
$\nabla^{g}\omega_\alpha=0$, for each
$\alpha\in\left\{1,\ldots,k\right\}$. We observe, firstly, that
$\nabla^g J_\alpha=0$, for each
$\alpha\in\left\{1,\ldots,k\right\}$. This is a consequence of the
definition of $J_\alpha$, of the fact that the leaves of the
foliation defined by $L_\alpha\oplus Q$ are K\"{a}hlerian manifolds,
and of the above properties that $\nabla^g L_\alpha\subset L_\alpha$
and $\nabla^g Q\subset Q$. Now we can prove that
$(\nabla^g_{V}\omega_\alpha)(W,W')=0$, for any
$V,W,W'\in\Gamma\left(TM\right)$. This equality holds immediately
for $W,W'\in\Gamma\left(L\right)$ and for
$W,W'\in\Gamma\left(Q\right)$ because $L$ and $Q$ are preserved by
$\nabla^g$. So it remains to show that $(\nabla^g_{V}
\omega_\alpha)\left(X,Y\right)=0$, for any
$X\in\Gamma\left(L\right)$ and $Y\in\Gamma\left(Q\right)$. In fact,
since $\nabla^g J_\alpha=0$ and $\nabla^g g=0$,
\begin{align*}
\left(\nabla^g_{V}\omega_\alpha\right)\left(X,Y\right)&=V\left(g\left(X,J_\alpha Y\right)\right)-g\left(\nabla^g_{V}X,J_\alpha Y\right)-g\left(X,J_\alpha \nabla^g_{V}Y\right)\\
&=V\left(g\left(X,J_\alpha Y\right)\right)-g\left(\nabla^g_{V}X,J_\alpha Y\right)-g\left(X,\nabla^g_{V}J_\alpha Y\right)\\
&=\left(\nabla^g_{V}g\right)\left(X,J_\alpha Y\right)=0.
\end{align*}
This concludes the proof.
\end{proof}

\section{Applications}

In this section we will examine some consequences of Theorem
\ref{connection}.  It can be useful to find the connection defined
in Theorem \ref{connection} in Darboux coordinates
$\{x_1,\ldots,x_n,$ $y_{1},\ldots,y_{kn}\}$ according to Theorem
\ref{darboux}. \ There exist functions $t^{\alpha j}_i$ such that
$Q=\textrm{span}\left\{X_1,\ldots,X_n\right\}$, where
$X_i:=\frac{\displaystyle \partial}{\displaystyle \partial
x_i}-\sum_{\alpha=1}^{k}\sum_{j=1}^{n}t_i^{\alpha
j}\frac{\displaystyle \partial}{\displaystyle \partial
y_{(\alpha-1)n+j}}$. We put $Y_{\alpha i}:=\frac{\displaystyle
\partial}{\displaystyle \partial y_{(\alpha-1)n+i}}$. Then by a straightforward
computation we have that
\begin{gather*}
\nabla_{Y_{\alpha i}}Y_{\beta j}=0, \ \
\nabla_{Y_{\alpha i}}X_j=0,\\
\nabla_{X_i}Y_{\alpha
j}=\sum_{\beta=1}^{k}\sum_{h=1}^{n}\frac{\displaystyle \partial
t_i^{\beta h}}{\displaystyle \partial y_{(\alpha-1)n+j}}Y_{\beta h},
\ \nabla_{X_i}X_j=-\sum\limits_{h=1}^{n}\frac{\displaystyle \partial
t^{\alpha j}_i}{\displaystyle \partial y_{(\alpha-1)n+h}}X_h,
\end{gather*}
where the functions $t^{\alpha j}_i$ satisfy the conditions
$\frac{\displaystyle \partial t^{\alpha j}_i}{\displaystyle \partial
y_{(\beta-1)n+h}}=0$ for $\alpha\neq\beta$, and $\frac{\displaystyle
\partial t^{1j}_i}{\displaystyle \partial y_{h}}=\cdots=\frac{\displaystyle \partial t^{kj}_i}{\displaystyle \partial
    y_{(k-1)n+h}}$, for all $i,j,h\in\left\{1,\ldots,n\right\}$. Moreover, the curvature
is given by
\begin{gather}
R_{Y_{\alpha i},Y_{\beta j}}=0, \ \ R_{X_i,X_j}=0, \label{significatocurvatura1} \\
R_{Y_{\alpha i},X_j}Y_{\beta
h}=\sum_{\gamma=1}^{k}\sum_{l=1}^{n}\frac{\displaystyle \partial^2
t^{\gamma_{l}}_j}{\displaystyle \partial y_{(\alpha-1)n+i}\partial
y_{(\beta-1)n+h}}Y_{\gamma l},\label{significatocurvatura2}\\
R_{Y_{\alpha i},X_j}X_k=-\sum_{l=1}^{n}\frac{\displaystyle
\partial^2 t^{\alpha_{k}}_i}{\displaystyle \partial y_{(\alpha-1)n+i}\partial
y_{(\alpha-1)n+l}}X_l.\label{significatocurvatura3}
\end{gather}
Then we have that the curvature $2$-form of $\nabla$ has the
following very simple expression\footnote{Throughout all this work,
if no confusion is feared, we identify forms on $M$ with their lifts
to principal bundle of linear frames $LM$.}:
\begin{equation*}
\Omega=\sum\Omega_{\alpha_i j}dx_{i}\wedge dy_{(\alpha-1)n+j}.
\end{equation*}
from which it follows that $\Omega^h$  vanishes for $h>n$. Thus if
$f\in I^h\left(G\right)$ is an $\textrm{ad}\left(G\right)$-invariant
polynomial of degree $h$, where $G=\emph{{Sp}}(k,n;\mathbb R)$, we
have that $f\left(\Omega\right)=0$ for $h=\deg\left(f\right)>n$.
This proves the following result.

\begin{proposition}
Under the assumptions of Theorem \ref{connection}, we have that
$\emph{Pont}^j\left(TM\right)=0$ for all $j>2n$, where
$\emph{Pont}(TM)$ denotes the Pontryagin algebra of the bundle $TM$.
\end{proposition}

Another strong consequence of Theorem \ref{connection} is the
existence of an Ehresmann connection. We recall the concept of
\emph{Ehresmann connection} for foliations. Let
$\left(M,\mathcal{F}\right)$ be a foliated manifold and $D$ a
distribution on $M$ which is supplementary to the tangent bundle $L$
of the foliation $\mathcal{F}$ at every point. A \emph{horizontal
curve} is a piecewise smooth curve
$\beta:\left[0,b\right]\longrightarrow M$, $b\in\mathbb{R}$, such
that $\beta'\left(t\right)\in D_{\beta\left(t\right)}$ for all
$t\in\left[0,b\right]$. A \emph{vertical curve} is a piecewise
smooth curve $\alpha:\left[0,a\right]\longrightarrow M$,
$a\in\mathbb{R}$, such that $\alpha'\left(t\right)\in
L_{\alpha\left(t\right)}$ for all $t\in\left[0,a\right]$, i.e. which
lies entirely in one leaf of $\mathcal{F}$. A \emph{rectangle} is a
piecewise smooth map
$\sigma:\left[0,a\right]\times\left[0,b\right]\longrightarrow M$
such that for every fixed $s\in\left[0,b\right]$ the curve
$\sigma_{s}:=\sigma|_{\left[0,a\right]\times \left\{s\right\}}$ is
vertical and for every fixed $t\in\left[0,a\right]$ the curve
$\sigma^t:=\sigma|_{\left\{t\right\}\times\left[0,b\right]}$ is
horizontal. The curves $\sigma_0=\sigma\left(\cdot,0\right)$,
$\sigma_b=\sigma\left(\cdot,b\right)$,
$\sigma^0=\sigma\left(0,\cdot\right)$ and
$\sigma^a=\sigma\left(a,\cdot\right)$ are called, respectively, the
\emph{initial  vertical  edge}, the \emph{final  vertical edge}, the
\emph{initial  horizontal  edge} and the \emph{final horizontal
edge} of $\sigma$. We say that the distribution  $D$ is  an
\emph{Ehresmann connection} for $\mathcal{F}$  if for every vertical
curve $\alpha$ and horizontal curve $\beta$ with the same initial
point, there exists a rectangle whose initial edges are $\alpha$ and
$\beta$ (cf. \cite{blumenthal1}). This rectangle is unique and is
called the \emph{rectangle associated to $\alpha$ and $\beta$}. It
is known (\cite{blumenthal0}) that every totally geodesic foliation
of a complete Riemannian manifold admits an Ehresmann connection,
namely the distribution orthogonal to the leaves of the foliation.
Furthermore, by the duality Riemannian -- totally geodesic, the
orthogonal bundle to a Riemannian foliation is also an Ehresmann
connection for this foliation.

Recall that given a foliated manifold $\left(M,\mathcal{F}\right)$
and a supplementary subbundle $D$ to $T\mathcal{F}$ (not necessarily
an Ehresmann connection), any horizontal curve
$\tau:\left[0,1\right]\longrightarrow M$ defines a family of
diffeomorphisms $\left\{\varphi_t:V_0\longrightarrow
V_t\right\}_{t\in\left[0,1\right]}$ such that
\begin{enumerate}
    \item each $V_t$ is a neighborhood of $\tau\left(t\right)$ in
    the leaf of $\mathcal{F}$ through $\tau\left(t\right)$, for all
    $t\in\left[0,1\right]$,
    \item
    $\varphi_t\left(\tau\left(0\right)\right)=\tau\left(t\right)$
    for all $t\in\left[0,1\right]$,
    \item for any fixed $p\in V_0$ the curve $t\mapsto
    \varphi_t\left(x\right)$ is horizontal,
    \item $\varphi_0:V_0\longrightarrow V_0$ is the identity map.
\end{enumerate}
This family of diffeomorphisms is called an \emph{element of
holonomy} along $\tau$ (\cite{blumenthal1}). It is shown in
(\cite{johnson}) and in (\cite{blumenthal0}) that an element of
holonomy along $\tau$ exists and is unique, in the sense that any
two elements of holonomy must agree on some neighborhood of
$\tau\left(0\right)$ in the leaf through $\tau\left(0\right)$. When
the leaves of $\mathcal F$ have a geometric structure -- such as a
Riemannian metric or a linear connection -- we say that \emph{$D$
preserves the geometry of the leaves} if the element of holonomy
along any horizontal curve is a local isomorphism of the particular
geometric structure.

Using the canonical connection which we have defined in $\S$
\ref{canonicalconnection} we prove now the following result.

\begin{theorem}\label{ehresmanntheorem}
Let $\left(M,\omega_\alpha,\mathcal F\right)$,
$\alpha\in\left\{1,\ldots,k\right\}$, be a compact connected
$k$-symplectic manifold and let $Q$ be an integrable distribution
transversal to $\mathcal F$ and satisfying the assumptions of
Theorem \ref{connection}. If the leaves of $\mathcal F$ are complete
affine manifolds, then the distribution $Q$ is an Ehresmann
connection for $\mathcal F$. Furthermore, if the canonical
connection $\nabla$ on $M$ induced by $Q$ is everywhere flat, then
the Ehresmann connection $Q$ preserves $\nabla$.
\end{theorem}
\begin{proof}
Let $\alpha:[0,a]\longrightarrow M$ and $\beta:[0,b]\longrightarrow
M$ be, respectively, a vertical and a horizontal curve such that
$\alpha(0)=x=\beta(0)$. We need to show that there exists a full
rectangle $\sigma:[0,a]\times[0,b]\longrightarrow M$ whose initial
vertical and horizontal edges are just $\alpha$ and $\beta$,
respectively. First we will show it under the further assumption
that $\alpha$ is a geodesic (with respect to the connection
$\nabla$). Fix an $s\in[0,b]$. We transport by parallelism the
vector $\alpha'(0)$ along the curve $\beta$, obtaining a vector
$v_s\in T_{\beta(s)}M$ which is in turn  tangent to $\mathcal F$
since the $\nabla$-parallel transport preserves the foliation
$\mathcal F$ (note also that the vector $v_s$ does not depend on the
curve because the curvature vanishes identically). Let $\tau_s$ be
the geodesic determined by the initial conditions
$\tau_s(0)=\beta(s)$ and $\tau'_s(0)=v_s$. Since the foliation
$\mathcal F$ is totally geodesic (with respect to $\nabla$),
$\tau_s$ is a curve lying on the leaf ${\mathcal L}_s$ of $\mathcal
F$ passing for $\beta(s)$, and the assumption on the completeness of
${\mathcal L}_s$ implies that we can extend $\tau_s$ for all the
values of the parameter $t$. In this way we obtain a map
$\sigma:[0,a]\times[0,b]\longrightarrow M$, defined by
$\sigma(t,s):=\tau_s (t)$, and it is easy to show that it is just
the rectangle we are looking for.  Now we have to prove the theorem
dropping the assumption that the curve $\alpha$ is a geodesic.
Because $M$ is  compact and the leaves of $\mathcal F$ are complete
affine manifolds with respect to $\nabla$, we find $\epsilon > 0$
such that for any $x \in M$, the $\epsilon$-ball $B(x, \epsilon)$ is
convex. As the leaves are totally geodesic, the $\epsilon$-balls
$B_{\mathcal L}(x, \epsilon)$ in any leaf ${\mathcal L}$ coincide
with the corresponding connected components of $B(x, \epsilon) \cap
{\mathcal L}$. Therefore, for any $x \in M$, there exists $\epsilon
> 0$ such that the $\epsilon$-balls $B_{\mathcal L}(x, \epsilon)$
are convex. Suppose now that $\alpha:\left[0,a\right]\longrightarrow
M$ is a vertical curve contained in
$B_{\mathcal{L}}\left(x,\epsilon\right)$, with
$x=\alpha\left(0\right)$. Let $\alpha_t$ denote the geodesic on
$\mathcal{L}$ joining $x$ with $\alpha\left(t\right)$, for any fixed
$t\in\left[0,a\right]$. Then we define
\begin{equation*}
\sigma\left(t,s\right):=\sigma_{\alpha_t,
\beta|_{\left[0,s\right]}}\left(t,s\right),
\end{equation*}
for any $\left(t,s\right)\in\left[0,a\right]\times\left[0,b\right]$,
where $\sigma_{\alpha_t, \beta|_{\left[0,s\right]}}$ denotes the
rectangle associated to the curves $\alpha_t$ and
$\beta|_{\left[0,s\right]}$. By the first part of the proof,
$\sigma$ is just the rectangle whose initial edges are $\alpha$ and
$\beta$. Finally, if $\alpha$ is \emph{any} leaf curve on $M$, not
necessarily contained in $B_{\mathcal{L}}\left(x,\epsilon\right)$,
then we can always find a partition of $\left[0,a\right]$, say
$0=t_0<t_1<\cdots<t_m=a$, with the property that, for any
$i\in\left\{0,\ldots,m-1\right\}$, $\alpha\left(t_i\right),
\alpha\left(t_{i+1}\right)\in
B\left(\alpha\left(t_i\right),\epsilon\right)$. Let $\sigma_{(0)}$
be the rectangle corresponding to $\alpha|_{\left[0,t_1\right]}$ and
$\beta$. The curve
$\beta_1:=\sigma_{\left(0\right)}|_{\left\{t_1\right\}\times\left[0,b\right]}$
is horizontal and $\beta_1\left(0\right)=\alpha\left(t_1\right)$, so
we can find a rectangle $\sigma_{(1)}$ whose edges are
$\alpha_{|_{\left[t_1,t_2\right]}}$ and $\beta_1$. After $m$ steps
we have $m$ rectangles
$\sigma_{(0)},\sigma_{(1)},\ldots,\sigma_{(m-1)}$ and we can define
$\sigma:=\sigma_{(0)}\cup\sigma_{(1)}\cup\cdots\cup\sigma_{(m-1)}$
obtaining the rectangle whose initial edges are $\alpha$ and
$\beta$. The last part of the statement follows directly from
(\cite[Proposition 5.3]{blumenthal1}).
\end{proof}

The existence of an Ehresmann connection implies strong consequences
for the foliation. Many of them have been studied in
(\cite{blumenthal1}), from which we have the following results.

\begin{corollary}
Let $(M,\omega_\alpha,{\mathcal F})$,
$\alpha\in\left\{1,\ldots,k\right\}$, be a $k$-symplectic manifold
satisfying the assumptions of Theorem \ref{ehresmanntheorem}. Then
the following statements hold:
\begin{enumerate}
  \item[(a)] Any two leaves of $\mathcal F$ can be joined by a
  horizontal curve.
  \item[(b)] The universal covers of any two leaves of $\mathcal{F}$ are
isomorphic.
  \item[(c)] The universal cover $\tilde{M}$ of $M$ is
topologically a product
$\tilde{\mathcal{L}}\times\tilde{\mathcal{Q}}$, where
$\tilde{\mathcal{L}}$ is the universal cover of the leaves of
$\mathcal{F}$ and $\tilde{\mathcal{Q}}$ the universal cover of the
leaves of the foliation integral to $Q$.
\end{enumerate}
\end{corollary}

In general, to each leaf $\mathcal{L}$ of a foliation admitting an
Ehresmann connection $D$, it is attached a group
$H_{D}\left({\mathcal{L}},x\right)$, $x\in\mathcal{L}$, defined as
follows (\cite{blumenthal1}). Let $\Omega_x$ be the set of all
horizontal curves $\beta:\left[0,1\right]\longrightarrow M$ with
starting point $x$. Then there is an action of the fundamental group
$\pi_1\left({\mathcal{L}},x\right)$ of $\mathcal{L}$ on $\Omega_x$
given in the following way: for any
$\delta=\left[\tau\right]\in\pi_1\left({\mathcal{L}},x\right)$ and
for any $\beta\in\Omega_x$, $\tau\cdot\beta$ is the final horizontal
edge of the rectangle corresponding to $\tau$ and $\beta$. It can be
proved that this definition does not depend on the vertical loop
$\tau$ in $x$ representing $\delta$. Let
$K_{D}\left({\mathcal{L}},x\right)=\left\{\delta\in\pi_1\left({\mathcal{L}},x\right):
\tau\cdot\beta=\beta \textrm{ for all } \beta\in\Omega_x\right\}$.
Then $K_{D}\left({\mathcal{L}},x\right)$ is a normal subgroup of
$\pi_1\left({\mathcal{L}},x\right)$ and we define
\begin{equation*}
H_{D}\left({\mathcal{L}},x\right):=\pi_1\left({\mathcal{L}},x\right)/K_{D}\left({\mathcal{L}},x\right).
\end{equation*}
It is known that $H_{D}\left({\mathcal{L}},x\right)$ does not depend
on the Ehresmann connection $D$, thus it is an invariant of the
foliation. Then we have the following result.

\begin{corollary}
Let $(M,\omega_\alpha,\mathcal F)$,
$\alpha\in\left\{1,\ldots,k\right\}$, be a $k$-symplectic manifold
satisfying Theorem \ref{ehresmanntheorem}. If $\mathcal{F}$ has a
compact leaf ${\mathcal{L}}_0$ with finite
$H_{D}\left({{\mathcal{L}}_0},x_0\right)$, then every leaf
$\mathcal{L}$ of $\mathcal{F}$ is compact with finite
$H_{D}\left({\mathcal{L}},x\right)$.
\end{corollary}
\begin{proof}
It is a direct consequence of (\cite[Theorem 1]{blumenthal2}).
\end{proof}

Another consequence of Theorem \ref{ehresmanntheorem} is the
following result.

\begin{corollary}\label{penultimo}
Let $(M,\omega_\alpha,\mathcal F)$,
$\alpha\in\left\{1,\ldots,k\right\}$, be a $k$-symplectic manifold
satisfying the assumptions of Theorem \ref{ehresmanntheorem}. Then
$\mathcal{F}$ has no vanishing cycles. Moreover, the homotopy
groupoid of $\mathcal{F}$ is a Hausdorff manifold.
\end{corollary}
\begin{proof}
The assertions follow from (\cite[Theorem 2]{wolak3}) and
(\cite[Corollary 2]{wolak3}).
\end{proof}

Now we study more deeply $k$-symplectic manifolds whose canonical
connections are flat. From
\eqref{significatocurvatura1}--\eqref{significatocurvatura3} it
follows that the geometric interpretation of the flatness of
$\nabla$ is that the functions $t^{\alpha_j}_i$ are leaf-wise
affine. Usually this condition is expressed saying that $Q$ is an
\emph{affine transversal distribution} for $\mathcal F$ (see, for
instance, \cite{vaisman1}, \cite{vaisman2}). In the following
theorem we give a normal form for flat $k$-symplectic manifolds:

\begin{theorem}\label{local}
Let $(M,\omega_\alpha,\mathcal F)$,
$\alpha\in\left\{1,\ldots,k\right\}$, be a $k$-symplectic manifold
and $Q$ a distribution satisfying the assumptions of Theorem
\ref{connection}. If the corresponding canonical connection $\nabla$
is flat, then there exist local coordinates
$\left\{x_1,\ldots,x_n,y_{1},\ldots,y_{kn}\right\}$ with respect to
which each $2$-form $\omega_\alpha$ is given by
\begin{equation}\label{darbouxforms}
\omega_\alpha=\sum_{i=1}^{n}dx_{i}\wedge dy_{(\alpha-1)n+i},
\end{equation}
$\mathcal F$ is described by the equations
$\left\{x_1=const.,\ldots,x_n=const.\right\}$ and $Q$ is spanned by
$\frac{\displaystyle \partial}{\displaystyle \partial
    x_{1}},\ldots,\frac{\displaystyle \partial}{\displaystyle \partial
    x_{n}}$.
\end{theorem}
\begin{proof}
Let $x\in M$ be a point and $U\subset M$ a chart containing $x$. One
can consider an adapted basis $\{e_{1},\ldots,e_{n(k+1)}\}$ of
$T_{x}M$ such that, for each $\alpha\in\left\{1,\ldots,k\right\}$,
$\{e_{(\alpha-1)n+1},\ldots,$ \ $e_{\alpha n}\}$ is a basis of
$L_{\alpha_{x}}$, $\{e_{kn+1},\ldots,e_{n(k+1)}\}$ is a basis of
$Q_{x}$, and
\begin{gather}\label{relazioni}
\omega_\alpha\left(e_{(\beta-1)n+i},e_{(\gamma-1)n+j}\right)=\omega_\alpha\left(e_{kn+i},e_{kn+j}\right)=0,\\
\omega_\alpha\left(e_{(\beta-1)n+i},e_{kn+j}\right)=-\frac{1}{2}\delta_{\alpha\beta}\delta_{ij},\label{relazioni1}
\end{gather}
for all $\alpha,\beta,\gamma\in\left\{1,\ldots,k\right\}$,
$i,j\in\left\{1,\ldots,n\right\}$. For each
$l\in\left\{1,\ldots,n(k+1)\right\}$ we define a vector field $E_k$
on $U$ by the $\nabla$-parallel transport along curves. More
precisely, for any $y\in U$ we consider a curve
$\gamma:\left[0,1\right]\longrightarrow U$ such that
$\gamma\left(0\right)=x$, $\gamma\left(1\right)=y$ and define
$E_l\left(y\right):=\tau_\gamma\left(e_l\right)$,
$\tau_\gamma:T_{x}M\longrightarrow T_{y}M$ being the parallel
transport along $\gamma$. Note that $E_{l}\left(y\right)$ does not
depend on the curve joining $x$ and $y$, since $R\equiv 0$. Thus we
obtain $n(k+1)$ vector fields on $U$, $E_{1},\ldots,E_{n(k+1)}$ such
that, for each $\alpha\in\left\{1,\ldots,k\right\}$,
$i\in\left\{1,\ldots,n\right\}$,
$E_{(\alpha-1)n+i}\in\Gamma(L_\alpha)$ and $E_{kn+i}\in\Gamma(Q)$,
since the connection $\nabla$ preserves the subbundles $L_\alpha$
and $Q$. Moreover, by \eqref{relazioni}--\eqref{relazioni1} we have
for any $y\in U$ and
$\alpha,\beta,\gamma\in\left\{1,\ldots,k\right\}$,
$i,j\in\left\{1,\ldots,n\right\}$
\begin{gather}\label{relazioni2}
\omega_\alpha\left(E_{(\beta-1)n+i},E_{(\gamma-1)n+j}\right)=\omega_\alpha\left(E_{kn+i},E_{kn+j}\right)=0,\\
\omega_\alpha\left(E_{(\beta-1)n+i},E_{kn+j}\right)=-\frac{1}{2}\delta_{\alpha\beta}\delta_{ij}.\label{relazioni3}
\end{gather}
Indeed, for all $l,m\in\left\{1,\ldots,n(k+1)\right\}$,
\begin{equation*}
\frac{d}{dt}\omega_\alpha\left(E_{l}\left(\gamma\left(t\right)\right),E_{m}\left(\gamma\left(t\right)\right)\right)=\omega_{\alpha}\left(\nabla_{\gamma'}E_{l},E_{m}\right)+\omega_{\alpha}\left(E_{l},\nabla_{\gamma'}E_{m}\right)=0
\end{equation*}
because $\omega_\alpha$ is parallel with respect to $\nabla$. Thus
$\omega_{{\alpha}_{x}}\left(e_{l},e_{m}\right)=\omega_{{\alpha}_{y}}\left(E_{l}\left(y\right),E_{m}\left(y\right)\right)$,
for any $y\in U$. Note that, by construction, we have
$\nabla_{E_l}E_m=0$ for all $l,m\in\left\{1,\ldots,n(k+1)\right\}$.
>From this, Theorem \ref{connection} and Proposition \ref{torsion} it
follows that the vector fields $E_1, \ldots, E_{n(k+1)}$ commute
each other. Therefore there exist local coordinates
$\{x_{1},\ldots,x_{n},y_{1},\ldots,y_{kn}\}$,
$\alpha\in\left\{1,\ldots,k\right\}$, such that
$E_{(\alpha-1)n+i}=\frac{\displaystyle \partial}{\displaystyle
\partial y_{\alpha_{i}}}$ and $E_{kn+j}=\frac{\displaystyle \partial}{\displaystyle
\partial
x_{j}}$, for any $i,j\in\left\{1,\ldots,n\right\}$. Note that by
\eqref{relazioni2}--\eqref{relazioni3} we get that
$\omega_\alpha=\sum_{i=1}^{n}dx_{i}\wedge dy_{(\alpha-1)n+i}$. Thus,
with respect this coordinate system,
\begin{enumerate}
    \item[(i)] each $L_\alpha$ is spanned by $\frac{\displaystyle \partial}{\displaystyle \partial
    y_{(\alpha-1)n+1}},\ldots,\frac{\displaystyle \partial}{\displaystyle \partial
    y_{\alpha n}}$,
    \item[(ii)] $Q$ is spanned by $\frac{\displaystyle \partial}{\displaystyle \partial
    x_{1}},\ldots,\frac{\displaystyle \partial}{\displaystyle \partial
    x_{n}}$,
    \item[(iii)] the $k$-symplectic forms $\omega_\alpha$ are given by $\omega_\alpha=\sum_{i=1}^{n}dx_{i}\wedge dy_{(\alpha-1)n+i}$.
\end{enumerate}
This proves the assertion.
\end{proof}

\begin{remark}
Theorem \ref{local} should be compared with Theorem \ref{darboux}.
It should be remarked that according to Theorem \ref{darboux} there
always exist local coordinates
$\{x_{1},\ldots,x_{n},y_{1},\ldots,y_{kn}\}$ verifying
\eqref{darbouxforms} and  such that the foliation $\mathcal F$ is
locally given by the equations
$\{x_1=\textrm{const.},\ldots,x_n=\textrm{const.}\}$. On the other
hand, by the general theory of foliations there always exists local
coordinates $\{x'_{1},\ldots,x'_{n},y'_{1},\ldots,y'_{kn}\}$ with
respect to which the foliation defined by $Q$ is described by the
equations $\{y'_1=\textrm{const.},\ldots,y'_{kn}=\textrm{const.}\}$.
In general these two types of coordinate systems do not coincide.
Theorem \ref{local} just states that a sufficient condition for this
is expressed by the flatness of the canonical connection. Note that
this condition is also necessary, as easily it follows from
\eqref{significatocurvatura1}--\eqref{significatocurvatura3}.
\end{remark}

\section*{Acknowledgments}
The authors thank the referees for their useful remarks.


\small

\begin{thebibliography}{99}

\bibitem{awane1}A. Awane, \textit{$k$-symplectic structures},
J. Math. Phys. \textbf{33} (1992), 4046--4052.

\bibitem{awane3}A. Awane, \textit{$G$-spaces $k$-symplectic
homog\`{e}nes}, J. Geom. Phys. \textbf{13} (1994), 139-157.

\bibitem{awane2}A. Awane, \textit{Some affine properties of the $k$-symplectic
manifolds}, Beitr. Algebra Geom. \textbf{39} (1998), 75--83.

\bibitem{awane5}A. Awane, M. Goze, \textit{Pfaffian systems, $k$-symplectic
systems}, Kluwer Academic Publishers, 2000.

\bibitem{blumenthal0}R. Blumenthal, J. Hebda, \textit{De Rham decomposition theorems for foliated
manifolds}, Ann. Inst. Fourier \textbf{33} (1983), 183--198.

\bibitem{blumenthal1}R. Blumenthal, J. Hebda, \textit{Ehresmann connections for foliations}, Indiana Univ. Math. J. \textbf{33} n. 4 (1984), 597--611.

\bibitem{blumenthal2}R. Blumenthal, J. Hebda, \textit{Complementary distributions which preserve the leaf geometry and applications to totally geodesic
foliations}, Quart. J. Math. Oxford \textbf{35} (1984), 383--392.

\bibitem{bott}R. Bott, \textit{Lectures on characteristic classes
and foliations}, Lect. Notes Math. \textbf{279} (1972), 1--94.

\bibitem{boyom}N. B. Boyom, \textit{M´etriques K¨ahl´eriennes affinement plates de certaines vari´et´es symplectiques I}, Proc. London
Math. Soc. \textbf{66} (1993), 358–80.

\bibitem{haggar}G. B. Byrnes, H. W. Capel, F. A. Haggar, G. R. W.
Quispel, \textit{$k$-integrals and $k$-Lie symmetries in discrete
dynamical systems}, Physica A \textbf{233} (1996), 379--394.

\bibitem{etayo1}F. Etayo, R. Santamaria, \textit{The canonical connection of a bi-Lagrangian
manifold}, J. Phys. A.: Math. Gen. \textbf{34} (2001), 981--987.

\bibitem{gunter}C. G\"{u}nther, \textit{The polysymplectic Hamiltonian formalism in field theory and calculus of variations
I: The local case}, J. Differential Geom. \textbf{25} (1987), 23-53.

\bibitem{hess}H. Hess, \textit{Connections on symplectic
manifolds and geometric quantization}, Lect. Notes Math.
\textbf{836} (1980), 153--166.

\bibitem{johnson}D. Johnson, L. Whitt - \textit{Totally geodesic
foliations}, J. Differential Geom. \textbf{15} (1980), 225--235

\bibitem{kobayashi}S. Kobayashi, K. Nomizu, \textit{Foundations of differential geometry, Vol. II}, Interscience Publishers, 1969.

\bibitem{deleon}M. de Le\'{o}n, E. Merino, J. A. Oubi\~{n}a, P. R.
Rodrigues, M. A. Salgado, \textit{Hamiltonian systems on
$k$-cosymplectic manifolds}, J. Math. Phys. \textbf{39} (1998),
876--893.

\bibitem{deleon1}M. de Le\'{o}n, M. McLean, L. K. Norris, A. Rey Roca, M. Salgado, \textit{Geometric structures in field
theory}, arXiv.math-ph/0208036 v1 (2002).

\bibitem{merino}E. Merino, \textit{$k$-symplectic and $k$-cosymplectic geometries. Applications to classical field
theory}. Pubblicaciones del Departemento de Geometria y Topologia,
Universitad de Santiago de Compostela \textbf{87}, 1997.

\bibitem{puta}M. Puta, \textit{Some remarks on the $k$-symplectic
manifolds}, Tensor \textbf{47} (1988), 109--115.

\bibitem{rey}A. M. Rey, N. Rom\'{a}n-Roy, M. Salgado, \textit{G\"{u}nter's formalism ($k$-symplectic formalism) in classical field theory: Skinner-Rusk approach and the evolution
operator}, J. Math. Phys. \textbf{46} (2005).

\bibitem{rovenski}V. Rovenskij, \textit{Foliations on Riemannian manifolds and
submanifolds}, Birkh\"{a}user, 1998.

\bibitem{vaisman1}I. Vaisman, \textit{$d_f$-cohomology of Lagrangian
foliations}, Monash. Math. \textbf{106} (1998), 221--244.

\bibitem{vaisman2}I. Vaisman, \textit{Basic of Lagrangian
foliations}, Publ. Mat. \textbf{33} (1989), 559--575.

\bibitem{wolak1}R. Wolak, \textit{Ehresmann connections for Lagrangian
foliations}, J. Geom. Phys. \textbf{17} (1995), 310--320.

\bibitem{wolak3}R. Wolak, \textit{Graphs, Ehresmann connections and vanishing
cycles}, Proc. Conf. Differential geometry and applications
(1996), 345--352.

\bibitem{wolak2}R. Wolak, \textit{On leaves of Lagrangian
foliations}, Russian Mathematics (Izvestiya VUZ. Matematika)
\textbf{42} n. 6 (1998), 24--31.
\end{thebibliography}
\end{document}